\documentclass[preprint,12pt]{elsarticle}

\usepackage{amssymb}
\usepackage{algorithmicx,algpseudocode,booktabs}
\usepackage{subfig}

\journal{Applied Mathematics and Computation}

\begin{document}

\begin{frontmatter}



\title{Nonlinear system identification and control using \\state transition algorithm}


\author{Xiaojun Zhou$^{\dagger\ddagger}$}
\author{Chunhua Yang$^{\dagger}$\corref{ychh}}
\author{Weihua Gui$^{\dagger}$}
\cortext[ychh]{Corresponding author. Tel.: +86 731 88836876.\\
\indent E-mail addresses: tiezhongyu2010@gmail.com (Xiaojun Zhou), ychh@csu.edu.cn (Chunhua Yang), gwh@csu.edu.cn (Weihua Gui). }
\address{$^{\dagger}$School of Information Science and Engineering,
Central South University,
Changsha 410083, China \\
$^{\ddagger}$School of Science, Information Technology and Engineering, University of Ballarat, Victoria 3353, Australia}

\begin{abstract}
By transforming identification and control for nonlinear system into optimization problems, a novel optimization method named state transition algorithm (STA) is introduced to solve the problems. In the proposed STA, a solution to a optimization problem is considered as a state, and the updating of a solution equates to a state transition, which makes it easy to understand and convenient to implement. First, the STA is applied to identify the optimal parameters of the estimated system with previously known structure. With the accurate estimated model, an off-line PID controller is then designed optimally by using the STA as well. Experimental results have demonstrated the validity of the methodology, and comparisons to STA with other optimization algorithms have testified that STA is a promising alternative method for system identification and control due to its stronger search ability, faster convergence rate and more stable performance.
\end{abstract}


\end{frontmatter}


\section{Introduction}
The identification and control of nonlinear system have been widely studied in recent years \cite{fei,pli,chen,jli}. Before the design of a controller, it is necessary to achieve system identification. In general, the process of system identification can be decomposed into two steps: the selection of an appropriate identification model (system structure) and an estimation of the model's parameters, of which, the parameter estimation plays a relatively more important role since a specific class of models that can best describe the real system can usually be derived by mechanism analysis of industrial processes \cite{bi1}.\\
\indent As for techniques of parameter estimation, approaches such as least-squares method, instrumental variable method, correlative function method, and maximum-likelihood method are widely used \cite{bi2,bi3}. Especially for the least-squares method, it has been successfully utilized to identify the parameters in static and dynamic systems \cite{bi4}. However, most of these techniques have some fundamental issues, including their dependence on unrealistic assumptions such as unimodal performance and differentiability of the performance function, and they are easily getting trapped into local optimum, because these methods are in essence local search techniques based on gradient. For example, the least-squares method is only suitable for the model structure possessing some linear property. Once the model structure exhibits nonlinear performance, this approach often fails in finding a global optimum and becomes ineffective \cite{bi1,bi2,bi3,bi5}. \\
\indent Fortunately, the modern intelligent optimization algorithms, such as genetic algorithm (GA) \cite{bi6,bi7}, particle swarm optimization (PSO) \cite{bi8,bi9}, are global search techniques based not on gradient, and they have been successfully applied in various optimization problems even with multimodal property. As a matter of fact, some intelligent optimization algorithms have been utilized in the field of nonlinear system identification and control. In \cite{bi10}, estimation of bar parameters with binary-coded genetic algorithm was studied, and it was verified that the GAs can produce better results than most deterministic methods. Genetic algorithm based parameter identification of a hysteretic brushless exciter model was proposed in \cite{bi11}. In \cite{bi3,bi12}, real-coded genetic algorithms were applied for nonlinear system identification and controller tuning, and the simulation examples demonstrated the effectiveness of the GA based approaches. Then, in \cite{bi1,bi2,bi5}, parameter estimation and control of nonlinear system based on adaptive particle swarm optimization were presented, and examples confirmed the validity of the method. Further more, in \cite{bi13}, identification of Jiles-Atherton model parameters using particle swarm optimization, and in \cite{bi14}, parameters identification for PEM fuel-cell mechanism model based on effective informed adaptive particle swarm optimization were put forwarded subsequently. All of these indicate that intelligent optimization techniques are alternatives for traditional methods including gradient descent, quasi-Newton, and Nelde-Mead's simplex methods.\\
\indent Although GA and PSO are alternative approaches for the problem, they always encounter premature convergence and
their convergence rates are not so satisfactory when dealing with some complex or multimodal functions\cite{bi15,bi16}.
State transition algorithm (STA) is a novel optimization method based on the concept of state and state transition recently, which originates from the thought of state space representation and space transformation \cite{bi17,bi18}. In STA, four special transformation operators are designed, and they represent different search functions in space, which makes STA easy to understand and convenient to implement. For continuous function optimization problems, STA has exhibited comparable search ability compared with other intelligent optimization algorithms \cite{xzhou2012,xzhou2013}. \\
\indent In this paper, the STA is firstly introduced to identify the optimal parameters of nonlinear system. Then, we will discuss the off-line PID controller design by adopting STA according to the estimated model. The PID control is popular due to its ease of use, good stability and simple realization. The key issue for PID controller design is the accurate and efficient tuning of PID control gains: proportional gain $K_p$, integral gain $K_i$ and derivative gain $K_d$. For adjusting PID controller parameters efficiently, many methods were proposed. The Ziegler-Nichols method is an experimental one that is widely used; however, this method needs certain prior knowledge on a plant model \cite{bi19}. Once tuning the controller by Ziegler-Nichols method, a good but not optimum system response will be gained. On the other hand, many artificial intelligence techniques such as neural networks, fuzzy systems and neural-fuzzy logic have been widely applied to the appropriate tuning of PID controller gains \cite{bi20}. Besides these methods, modern intelligent optimization algorithms, such as GA and PSO, have also received much attention, and they are used to find the optimal parameters of PID controller \cite{bi1,bi3,bi12}.\\
\indent The goal of this paper is to introduce a novel method STA for both parameter estimation and control of nonlinear systems. In order to evaluate the performance of the STA, experiments are carried out to testify the validity of the proposed methodology, the results of which have confirmed that STA is an efficient method. Compared with other intelligent optimization algorithms, the simulation examples have demonstrated that the STA has superior features in terms of search ability, convergence rate and stability.



\section{Problems description}
To transform a specified problem into the standard form of optimization problem is called optimization modeling, which is the basis for parameter identification and system control. The standard optimization problems should consist of objective function and decision variables, while optimization algorithms are used to find a global optimal solution to the objective function restricted to some additional constraints.
\subsection{Identification of nonlinear system}
\begin{figure}[!htbp]
\centering
\includegraphics[width=8cm,height=6cm]{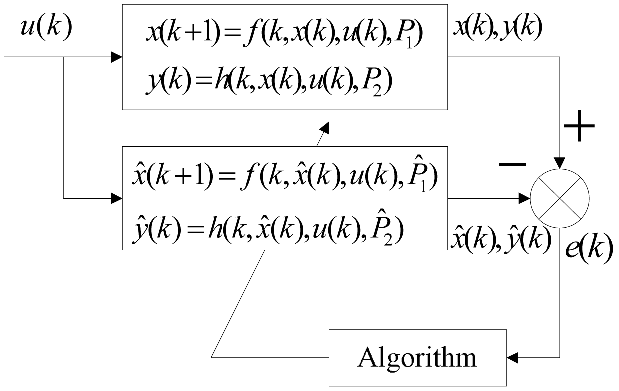}
\caption{The parameter estimation process}
\label{fig1}
\end{figure}

In this paper, the following class of discrete nonlinear systems described by the state space model is considered:
\begin{equation}
\begin{array}{ll}
x(k+1) & = f(k,x(k),u(k),P_1)\\
y(k)   & = h(k,x(k),u(k),P_2),
\end{array}
\end{equation}
where, $x \in \Re^n$ is the state vector, $u \in \Re$ is the input, $y \in \Re$ is the output, $P_1$ and $P_2$ are unknown parameter vectors that will be identified, and $f(\cdot)$ and $h(\cdot)$ are nonlinear functions. Without loss of generality, let $\theta = [\theta_1,\theta_2,\cdots,\theta_n]$ be a rearranging vector containing all parameters in $P_1$ and $P_2$ where $n$ represents the total number of unknown system parameters. Furthermore, the estimated system model can be described as:
\begin{equation}
\begin{array}{ll}
\hat{x}(k+1) & = f(k,\hat{x}(k),u(k),\hat{P_1})\\
\hat{y}(k)   & = h(k,\hat{x}(k),u(k),\hat{P_2}),
\end{array}
\end{equation}
where, $\hat{x} \in \Re^n$ and $\hat{y} \in \Re$ denote the state vector and the output of the model, $\hat{P_1}$ and $\hat{P_2}$ are the estimated parameter vectors, respectively. Accordingly, let $\hat{\theta} = [\hat{\theta}_1,\hat{\theta}_2,\cdots,\hat{\theta}_n]$ be the estimated rearranging vector.\\
\indent The basic thought of system identification is to compare the real system responses with the estimated system responses. Moreover, to accurately estimate the $\hat{\theta}$, some assumptions on the nonlinear systems are required:\\
(1) The system output must be available for measurement.\\
(2) System parameters must be connected with the system output.\\
\indent To deal with the problem of parameter estimation, a specified problem should be formulated as an optimization problem. In this study, the decision variables are the estimated parameter vector $\hat{\theta}$, while the objective function is chosen as the following mean squared errors(MSE):
\begin{equation}
\textrm{MSE} = \frac{1}{N}\sum_{k=1}^{N}e^2 = \frac{1}{N} \sum_{k=1}^{N} [X(k) - \hat{X}(k)]^2,
\end{equation}
where, $N$ is the length of sampling data, $X(k) = [x(k),y(k)]$ and $\hat{X}(k) = [\hat{x}(k),\hat{y}(k)]$ are real and estimated values at time $k$, respectively.\\
\indent It is obvious to find that the MSE is the function of variable vector $\hat{\theta}$, and then, the optimization problem will be solved by optimization algorithms which will minimize the MSE value so that the real nonlinear system is actually estimated. The block diagram of the nonlinear system parameter estimation is given in Fig.\ref{fig1}.
\subsection{Design of PID controller}
When the nonlinear system model is estimated, an off-line PID controller is then designed to guarantee the stability and other performances of the system. The reason why the PID controller is adopted is that it is the most widely used controller for application in industrial processes. The continuous form of a PID controller can be described as follows:
\begin{equation}
u(t) = K_p[e(t)+\frac{1}{T_i}\int_0^te(t)dt+T_d\frac{d}{dt}e(t)],
\end{equation}
where, $e(t)$ is the error signal between the desired and actual outputs, $u(t)$ is the control force, $K_p,T_i,T_d$ are the proportional gain, integral time constant and derivative time constant, respectively. By using the following approximations:
\begin{equation}
\begin{array}{ll}
\int_0^te(t)dt \approx T \sum_{j=0}^ke(j)\\
\frac{d}{dt}e(t) \approx \frac{e(k)-e(k-1)}{T},
\end{array}
\end{equation}
where, $T$ is the sampling period, then (5) can be rewritten as
\begin{equation}
u(k) = K_p\{e(k)+\frac{T}{T_i}\sum_{j=0}^ke(j)+\frac{T_d}{T}[e(k)-e(k-1)]\},
\end{equation}
which is called the place type, and in most case, the increment style as described following is more practical:
\begin{small}
\begin{equation}
\begin{array}{ll}
u(k) &= u(k-1)+ K_p\{[e(k)-e(k-1)]+\frac{T}{T_i}e(k)+\frac{T_d}{T}[e(k)-2e(k-1)+e(k-2)]\}\\
     &= u(k-1)+ K_p[e(k)-e(k-1)] + K_ie(k)+ K_d[e(k)-2e(k-1)+e(k-2)],
\end{array}
\end{equation}
\end{small}
where, $K_i$ and $K_d$ are integral gain and derivative gain, respectively.\\
\begin{figure}[!htbp]
\centering
\includegraphics[width=10cm,height=6cm]{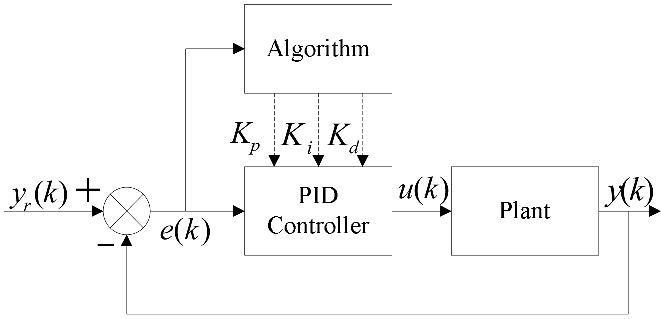}
\caption{The design of PID controller}
\label{fig2}
\end{figure}
\label{}
\indent The block diagram of the design of an off-line PID controller is illustrated in Fig.\ref{fig2}, where $y_r(k)$ is the reference output, and $y(k)$ is the system output at the sampling point. Optimization algorithms are used to adjust the PID controller parameters such as $K_p, K_i$ and $K_d$. In the same way, mean squared errors will be defined as the objective function
\begin{equation}
\textrm{MSE} = \frac{1}{N}\sum_{k=1}^{N}e^2 = \frac{1}{N}\sum_{k=1}^{N}[y_r(k) - y(k)]^2.
\end{equation}
\section{State transition algorithm}
Let's consider the following unconstrained optimization problem:
\begin{equation}
\min_{x \in \Re^{n}} f(x)
\end{equation}
where, $f:\Re^{n} \rightarrow \Re$ is a objective function. In a numerical way, the iterative method is adopted to solve the problem, the essence of which is to update the solution found so far. When thinking in a state and state transition way, a solution can be regarded as a state, and the updating of a solution can be considered as a state transition process. \\
\indent Based on the thought stated above, the form of state transition algorithm can be described as follows,
\begin{equation}
\left \{ \begin{array}{ll}
x_{k+1}= A_{k}x_{k} + B_{k}u_{k}\\
y_{k+1}= f(x_{k+1})
\end{array} \right.
\end{equation}
where, $x_{k}$ stands for a state, corresponding to a solution of a optimization problem; $A_{k}$ and
$B_{k}$ are state transition matrices, which are usually transformation operators;
$u_{k}$ is the function of variable $x_{k}$ and historical states; $f$ is the objective function or evaluation function.\\
\indent Using various types of space transformation for reference, four special
state transformation operators are designed to solve continuous function optimization problems.\\
(1) Rotation transformation
\begin{equation}
x_{k+1}=x_{k}+\alpha \frac{1}{n \|x_{k}\|_{2}} R_{r} x_{k},
\end{equation}
where, $x_{k}$ $\in$ $\Re^{n}$, $\alpha$ is a positive constant, called rotation factor;
$R_{r}$ $\in$ $\Re^{n\times n}$ is a random matrix with its elements belonging to the range of [-1, 1]
and $\|\cdot\|_{2}$ is 2-norm of a vector. It has proved that the rotation transformation
has the function of searching in a hypersphere \cite{bi17,bi18}.\\
(2) Translation transformation\\
\begin{equation}
x_{k+1} = x_{k}+  \beta  R_{t}  \frac{x_{k}-x_{k-1}}{\|x_{k}-x_{k-1}\|_{2}},
\end{equation}
where, $\beta$ is a positive constant, called translation factor; $R_{t}$ $\in \Re$ is a random variable
with its elements belonging to the range of [0,1]. It has illustrated
the translation transformation has the function of searching along a line from $x_{k-1}$ to $x_{k}$
 at the starting point $x_{k}$, with the maximum length of $\beta$ \cite{bi17,bi18}.\\
(3) Expansion transformation
\begin{equation}
x_{k+1} = x_{k}+  \gamma  R_{e}x_{k},
\end{equation}
where, $\gamma$ is a positive constant, called expansion factor; $R_{e} \in \Re^{n \times n}$ is a random diagonal
matrix with its elements obeying the Gaussian distribution. It has also stated the expansion transformation
has the function of expanding the elements in $x_{k}$ to the range of [-$\infty$, +$\infty$], searching in the whole space \cite{bi17,bi18}.\\
(4) Axesion transformation
\begin{equation}
x_{k+1} = x_{k}+  \delta  R_{a}  x_{k},
\end{equation}
where, $\delta$ is a positive constant, called axesion factor; $R_{a}$ $\in \Re^{n \times n}$ is a random diagonal matrix with its elements obeying the Gaussian distribution and only one random index has nonzero value. As illustrated in \cite{bi18}, the axesion transformation is aiming to search along the axes.\\
\indent When using these transformation operators into practice, an important parameter called search enforcement(SE) is introduced to describe the times of certain transformation.\\
\indent The main procedures of the version of state transition algorithm in \cite{bi18} can be outlined in the following pseudocode.
\begin{algorithmic}[1]
\Repeat
    \If{$\alpha < \alpha_{\min}$}
    \State {$\alpha \gets \alpha_{\max}$}
    \EndIf
    \State {Best $\gets$ expansion(funfcn,Best,SE,$\beta$,$\gamma$)}
    \State {Best $\gets$ rotation(funfcn,Best,SE,$\alpha$,$\beta$)}
    \State {Best $\gets$ axesion(funfcn,Best,SE,$\beta$,$\delta$)}
    \State {$\alpha \gets \frac{\alpha}{\textit{fc}}$}
\Until{the maximum iterations is met}
\end{algorithmic}
where, \textit{fc} is a constant coefficient used for lessening the $\alpha$, and the translation operator will only be performed when a better solution is obtained.
\label{}
\section{Experimental results and comparison}
For both the identification of nonlinear system and the design of an off-line PID controller, when the optimization algorithms are utilized, they help to minimize the mean squared errors (MSE). In other words, the MSE or the evaluation of the objective function will guide the search of the algorithms. Different from methods based on gradient, the termination criterion of intelligent optimization algorithms usually are not the precision of the gradient but a prespecified maximum number of iterations.\\
\indent For comparison, the maximum iterations, population size or search enforcement are the same, and they are fixed at 100 and 30, respectively. To be more specific, in PSO, $c_1 = c_2 =1$, and $w$ will decrease in a linear way from 0.9 to 0.4, as suggested in \cite{bi6}. In STA \cite{bi18}, the rotation factor $\alpha$ will decrease in an exponential way with base $f_c =2$ from $\alpha_{max} =1$ to $\alpha_{min} =1e^{-4}$, and translation factor $\beta$, expansion factor $\gamma$,   axesion factor $\delta$ are all constant at 1. For GA, we use the MATLAB genetic algorithm toolbox v1.2 from http://www.sheffield.ac.uk/acse/research/ecrg/getgat.html. In this paper, the following two instances are studied.\\
\\
\textbf{Example 1.} An unstable nonlinear system is described by
\begin{equation}
\begin{array}{ll}
x_1(k+1)= \theta_{1}x_1(k)x_2(k),x_1(0) = 1,\\
x_2(k+1)= \theta_{2}x_1^2(k)+u(k),x_2(0)=1,\\
y(k)= \theta_{3}x_2(k)-\theta_4x_1^2(k),
\end{array}
\end{equation}
where, $\theta_1,\cdots,\theta_4$ are to be estimated. The real parameters of the nonlinear system are assumed to be $\theta = [\theta_1,\theta_2,\theta_3,\theta_4]=[0.5,0.3,1.8,0.9]$.\\
\indent The relative variables used in optimization algorithms are given as follows\\
\indent \qquad \qquad \quad $\theta_1 \in [0,2], \theta_2 \in [0,2], \theta_3 \in [0,2], \theta_4 \in [0,2]$, $N=8. $\\
\indent Considering the randomness of the stochastic optimization algorithms, 30 independent trials are run. In the meanwhile, some statistics, such as \textit{best} (the minimum), \textit{mean}, \textit{worst} (the maximum), \textit{st.dev} (standard deviation), are used to evaluate the performance of the algorithms.\\
\begin{table}[!htbp]
\begin{center}
\caption{Best estimated parameters for Example 1}
\label{tab1}
\begin{tabular}{{ccccc}}
\hline
\toprule[1pt]
\textit{Algorithms} & $\theta_1$ & $\theta_2$ & $\theta_3$ & $\theta_4$\\
\hline
GA  & 0.4981  &  0.2995  &   1.7946 &    0.8946 \\
PSO & 0.5000  &  0.3000  &   1.8000 &    0.9000 \\
STA & 0.5000  &  0.3000  &   1.8000 &    0.9000 \\
\bottomrule[1pt]
\hline
\end{tabular}
\end{center}
\end{table}
\begin{table}[!htbp]
\begin{center}
\caption{Performance comparison for Example 1}
\label{tab2}
\begin{tabular}{{ccccc}}
\hline
\toprule[1pt]
\textit{Algorithms} & \textit{Best} & \textit{Mean} & \textit{Worst} & \textit{St.dev}\\
\hline
GA  & 9.6674E-07  &  1.2572E-04  &   4.3492E-04 &    1.4433E-04\\
PSO & 1.3938E-12  &  2.8000E-03  &   2.7700E-02 &    8.4000E-03\\
STA & 5.2364E-12  &  5.2367E-11  &   1.3729E-10 &    3.5120E-11\\
\bottomrule[1pt]
\hline
\end{tabular}
\end{center}
\end{table}

Table \ref{tab1} lists the best estimated parameters gained by GA, PSO and STA, from which, we can find that only PSO and STA can achieve the accurate parameters of the real system, and the results obtained by GA are a little far from the real parameters. Then, from Table \ref{tab2}, it indicates that STA is the most stable algorithm for the problem because the \textit{mean} and \textit{st.dev} of STA are the smallest. It can also be found the results gained by GA and PSO are not so satisfactory since the \textit{best} is deviated from the \textit{mean} seriously.

\begin{figure}[!htbp]
\centering
\subfloat[$\theta_1$]{\includegraphics[width=8cm,height=6cm]{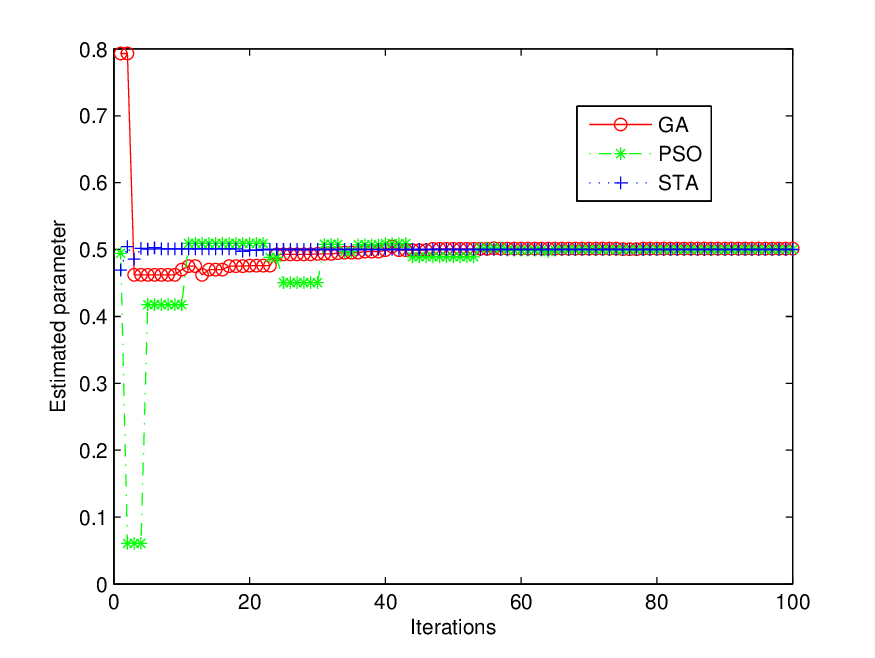}}
\subfloat[$\theta_2$]{\includegraphics[width=8cm,height=6cm]{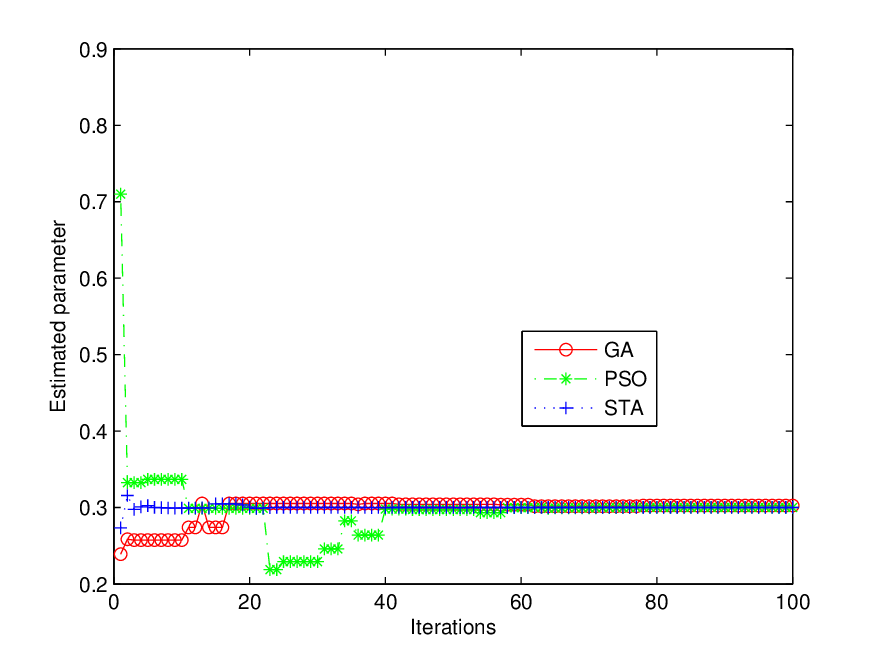}}\\
\subfloat[$\theta_3$]{\includegraphics[width=8cm,height=6cm]{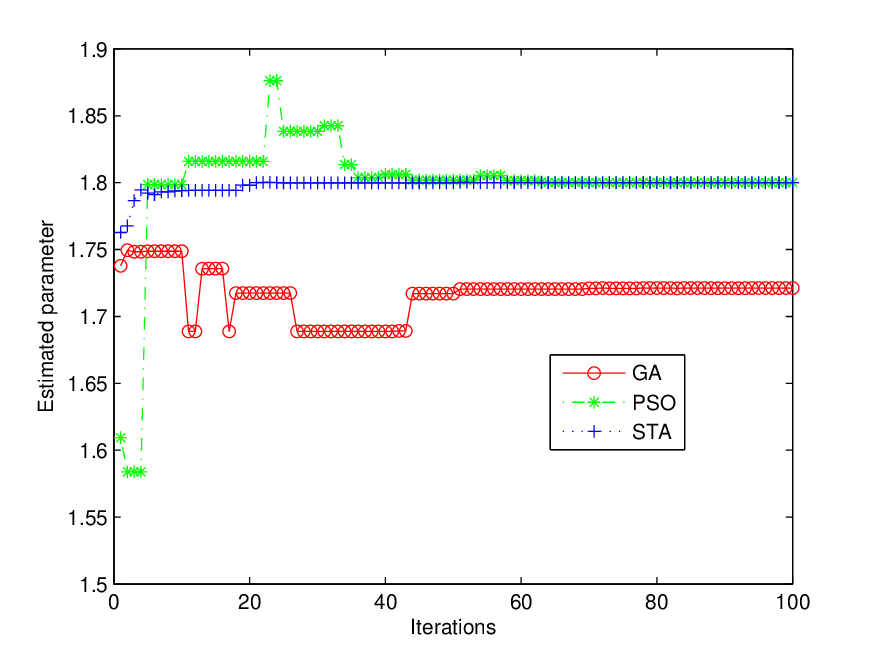}}
\subfloat[$\theta_4$]{\includegraphics[width=8cm,height=6cm]{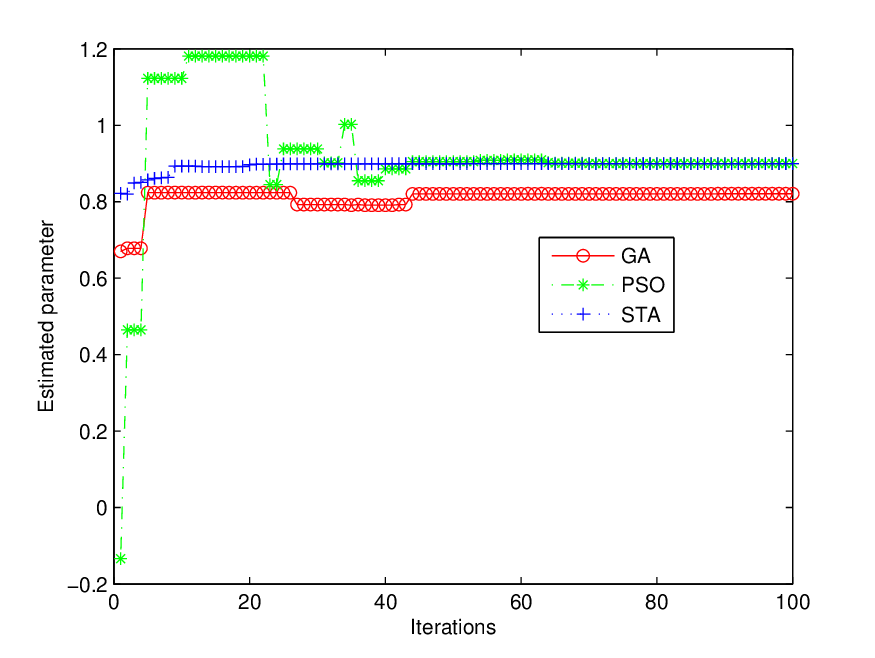}}
\caption{Trajectories of parameters $\theta_1,\theta_2,\theta_3$, and $\theta_4$ for Example 1.}
\label{fig3}
\end{figure}

Fig.\ref{fig3} illustrates the optimization processes of parameter estimation by using STA compared with other two algorithms in a middle run. It is easy to find that the convergence rate of STA are much faster than that of GA and PSO, with no more than 30 iterations, and the changes of parameters with STA are also steadier than other two algorithms.

\begin{figure}[htbp]
\centering
\subfloat[$K_p$]{\includegraphics[width=8cm,height=6cm]{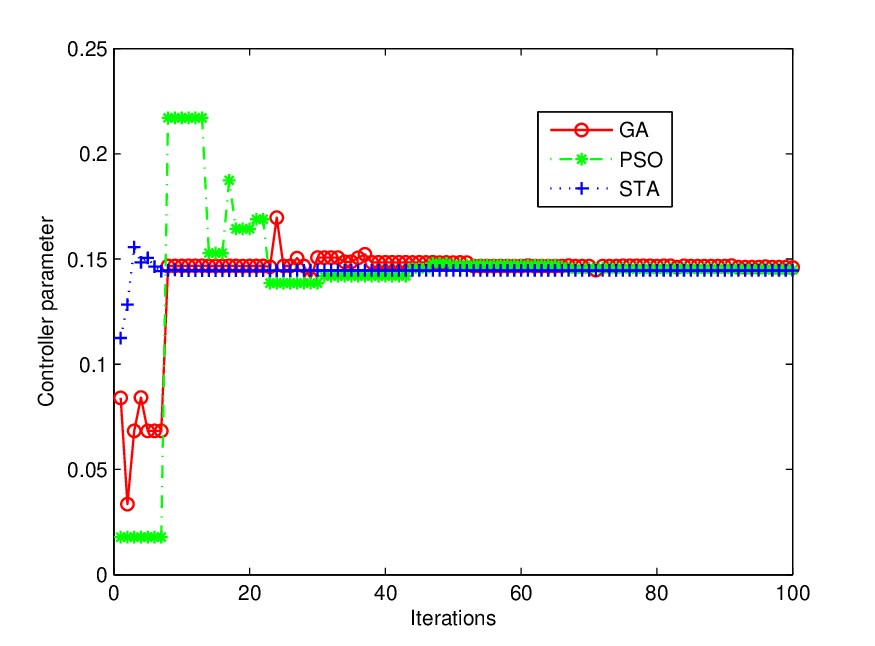}}
\subfloat[$K_i$]{\includegraphics[width=8cm,height=6cm]{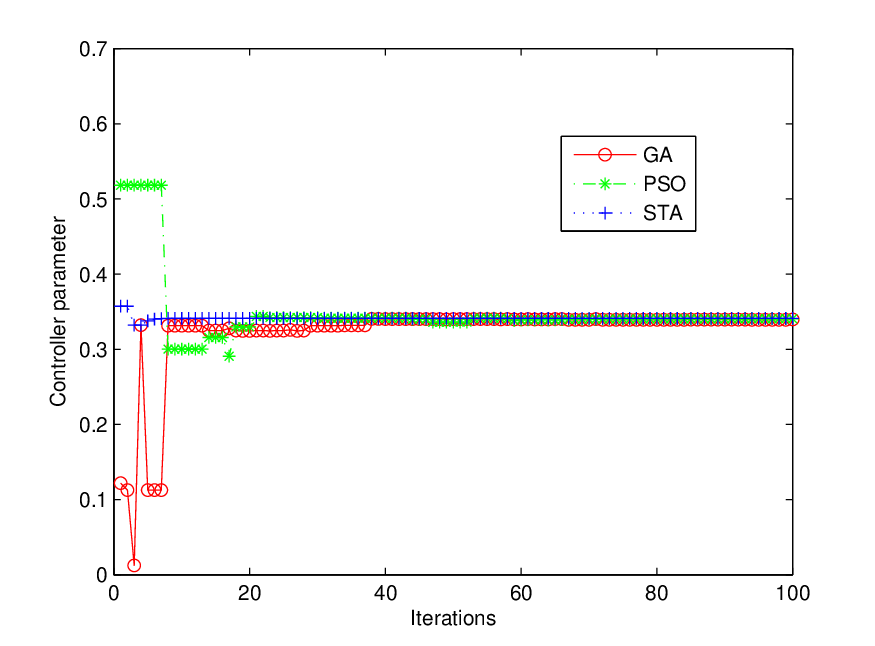}}\\
\subfloat[$K_d$]{\includegraphics[width=8cm,height=6cm]{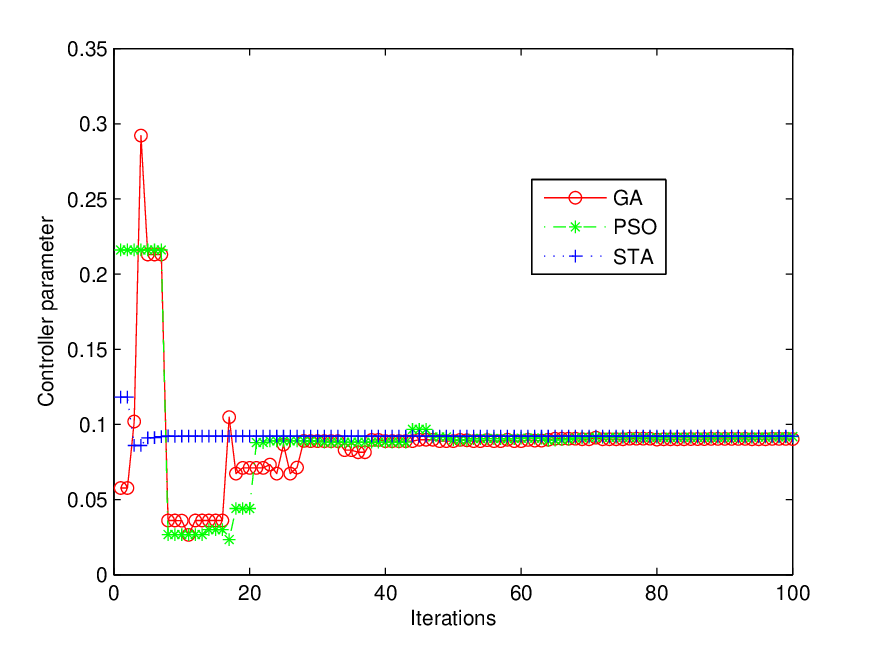}}
\subfloat[states and output]{\includegraphics[width=8cm,height=6cm]{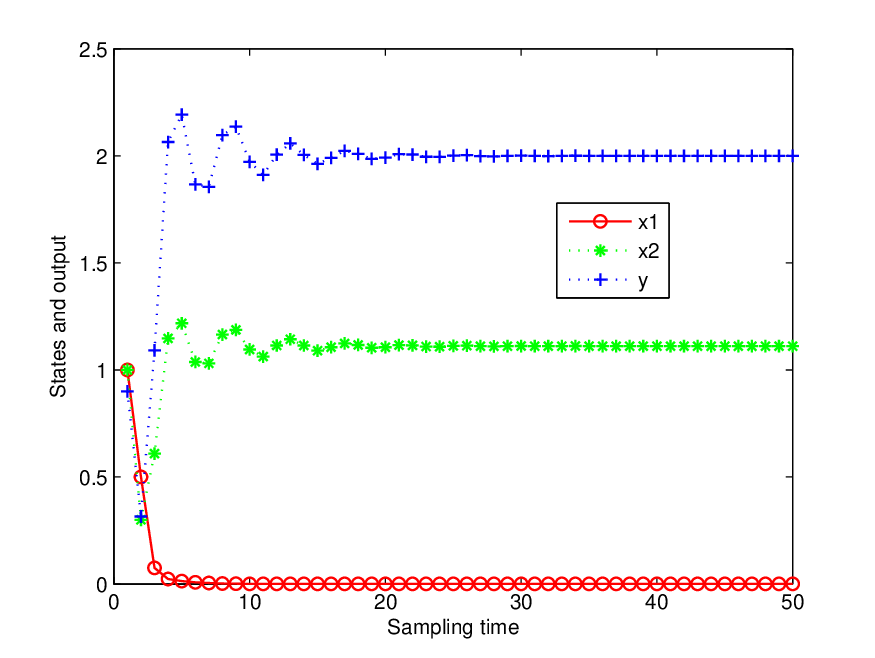}}\\
\caption{(a),(b),(c) depict the convergence of PID controller parameters $K_p,K_i$ and $K_d$ with GA, PSO and STA, respectively, and (d) shows the changes of states and output under the best parameters with STA for Example 1.}
\label{fig4}
\end{figure}

Then, with the estimated model, an off-line PID controller for this system is designed by using STA, too. As a optimizer, STA is to minimize the mean squared error between the plant output $y$ and the desired output $y_r$. In the experiment, relative variables are given by\\
\indent \qquad \qquad \quad $K_p \in [0,1], K_i \in [0,1], K_d \in [0,1]$, $y_r=2, N=50.$\\
\indent Also, in the same time, 30 independent trials are carried out, and some statistics are used to describe the performance of the algorithms. Table \ref{tab3} shows the parameters of PID controller under the best performance, and Table \ref{tab4} gives the detailed statistical results obtained by the three algorithms. Compared with GA and PSO, from Table \ref{tab4}, it is obvious to find that STA can find the minimum MSE with a higher probability, which shows that STA is much more appropriate for the problem.\\
\indent Fig.\ref{fig4} illustrates the convergence processes of PID controller parameters and the changes of states and output under best MSE using STA, which shows that the convergence rate is fast, and the changes of output indicate that the stability of nonlinear system is good under the proposed PID controller. To be more specific, for the nonlinear system, it can be found that $x_1$ is stable at 0, $x_2$ is stable at 1.1111, and $y$ is stable at 2 finally, which can be understood easily by analyzing the system.
\begin{table}[htbp]
\begin{center}
\caption{Best parameters of PID controller for Example 1}
\label{tab3}
\begin{tabular}{{cccc}}
\hline
\toprule[1pt]
\textit{Algorithms} & $K_p$ & $K_i$ & $K_d$\\
\hline
GA  & 0.1459  &  0.3398  &   0.0908\\
PSO & 0.1447  &  0.3407  &   0.0919\\
STA & 0.1445  &  0.3410  &   0.0922\\
\bottomrule[1pt]
\hline
\end{tabular}
\end{center}
\end{table}

\begin{table}[htbp]
\begin{center}
\caption{Performance comparison of PID controller for Example 1}
\label{tab4}
\begin{tabular}{{ccccc}}
\hline
\toprule[1pt]
\textit{Algorithms} & Best & Mean & Worst & St.dev\\
\hline
GA  & 0.1000 &  0.1002 &   0.1013 &    0.0017\\
PSO & 0.1000 &  0.1127 &   0.3089 &    0.1986\\
STA & 0.1000 &  0.1000 &   0.1000 &    1.41E-10\\
\bottomrule[1pt]
\hline
\end{tabular}
\end{center}
\end{table}

\textbf{Example 2.} Consider a first-order with time-delay system whose transfer function is given by
\begin{equation}\label{ex2}
G(s) = \frac{y(s)}{u(s)} = \frac{K}{Ts + 1} e^{-\tau s},
\end{equation}
where, $K$ is the steady-state gain, $T$ is the time constant, and $\tau$ is the time delay.
With the sampling time 0.01, the system (\ref{ex2}) can be easily changed to the discrete dynamic equation as follows:
\begin{equation}
\begin{array}{ll}
x(k+1)= (1 - \frac{1}{10T})x(k) + \frac{K}{10T}u(k-10\tau),\; x(0) = 0, \\
y(k) = x(k).
\end{array}
\end{equation}

It is assumed that the actual numerical values are $K = 10$, $T = 5$, and $\tau = 9$, respectively. Let consider $\theta = [\theta_1, \theta_2, \theta_3] = [K,T,\tau]$ as a vector of estimated parameters and set the control input $u(k) = 1$ in this study. The relative variables used in optimization algorithms are given as follows\\
\indent \qquad \qquad \quad $K \in [0,20], T \in [0,20], \tau \in [0,20]$, $N=350$.
\begin{table}[!htbp]
\begin{center}
\caption{Best estimated parameters for Example 2}
\label{tab5}
\begin{tabular}{{cccc}}
\hline
\toprule[1pt]
\textit{Algorithms} & $K$ & $T$ & $\tau$\\
\hline
GA  & 9.9710  &  4.9813  &   9.0283 \\
PSO & 9.9989  &  4.9993  &   9.0414 \\
STA & 10.0000 &  5.0000  &   9.0000 \\
\bottomrule[1pt]
\hline
\end{tabular}
\end{center}
\end{table}
\begin{table}[!htbp]
\begin{center}
\caption{Performance comparison for Example 2}
\label{tab6}
\begin{tabular}{{ccccc}}
\hline
\toprule[1pt]
\textit{Algorithms} & \textit{Best} & \textit{Mean} & \textit{Worst} & \textit{St.dev}\\
\hline
GA  & 2.0759E-07  &  0.0059       &   0.0422     &  0.0083\\
PSO & 2.2503E-10  &  0.0113       &   0.2412     &  0.0436\\
STA & 2.4004E-19  &  1.2816E-16   &   2.4510E-15 &  4.6216E-16\\
\bottomrule[1pt]
\hline
\end{tabular}
\end{center}
\end{table}

\begin{figure}[htbp]
\centering
\subfloat[$K$]{\includegraphics[width=8cm,height=6cm]{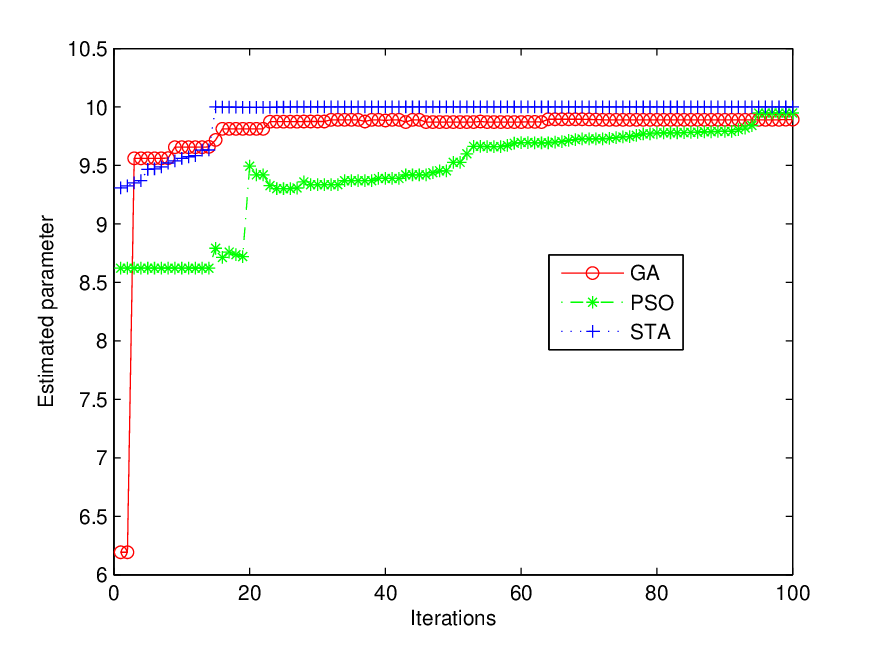}}
\subfloat[$T$]{\includegraphics[width=8cm,height=6cm]{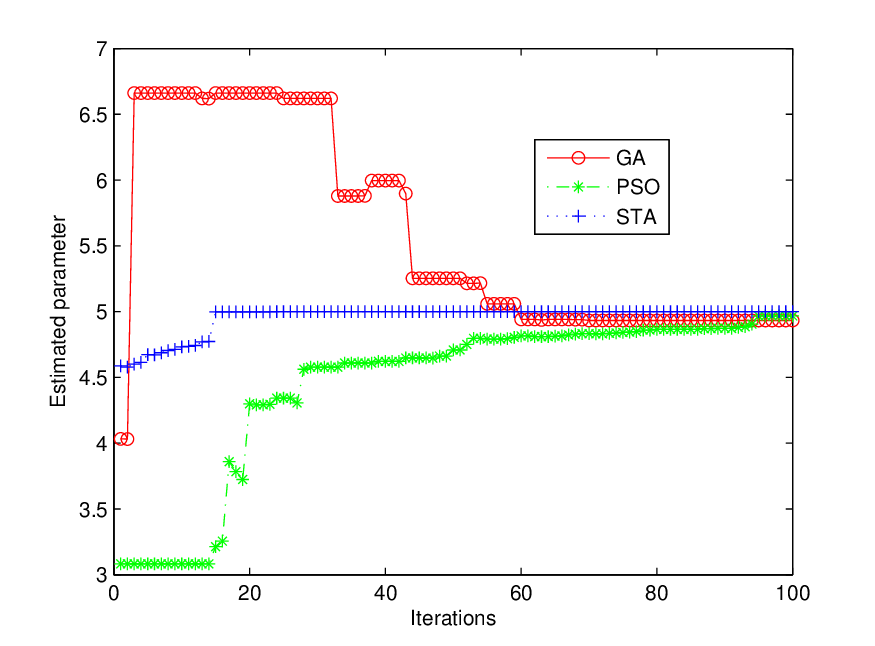}}\\
\subfloat[$\tau$]{\includegraphics[width=8cm,height=6cm]{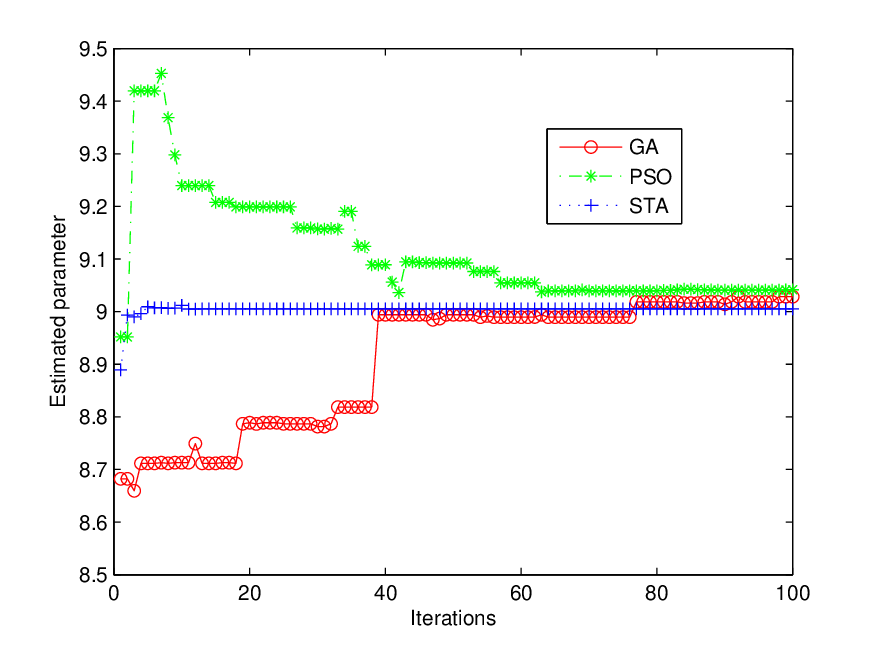}}
\subfloat[MSE]{\includegraphics[width=8cm,height=6cm]{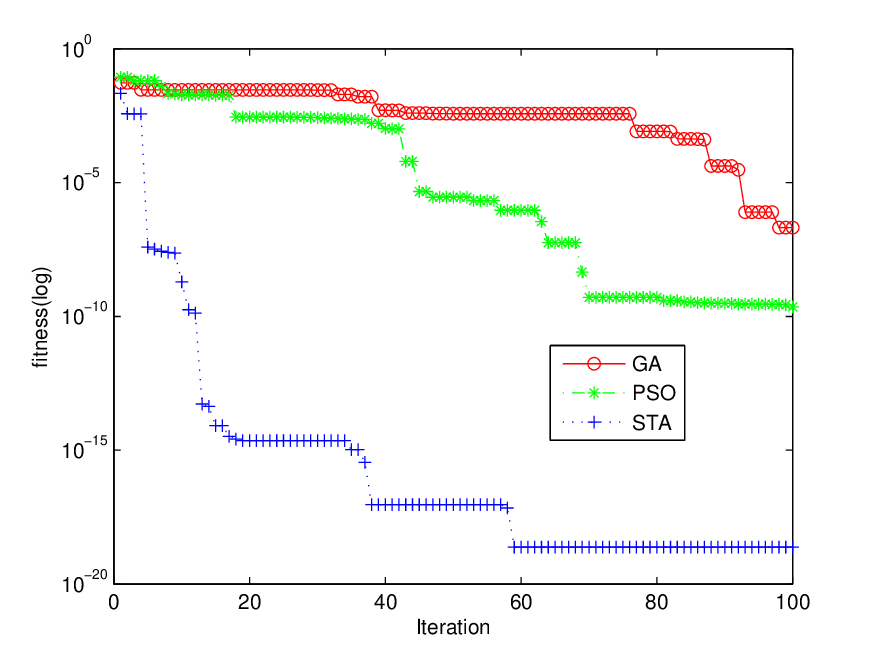}}
\caption{(a),(b),(c) show trajectories of parameters $K,T$, and $\tau$, and (d) shows the convergence of MSE for Example 2 under GA, PSO and STA, respectively.}
\label{fig5}
\end{figure}

Table \ref{tab5} shows the best estimated parameters found by GA, PSO and STA, respectively. Considering the
reformulated optimization problem is highly nonlinear and nonconvex, only STA can achieve the global solution of the real system with high precision, that is to say, STA has much stronger global search ability and much higher calculation accuracy than its competitors.  Anyway, both GA and PSO have also found the approximately global solution, as indicated in Table \ref{tab6}, since the \textit{best} of GA and PSO are approaching zero closely. However, the \textit{mean}, the \textit{worst} and the \textit{st.dev} are much far away from zeroes,
which indicates that GA and PSO are not stable for this identification problem. On the other hand, the \textit{worst} of STA is very small, which demonstrates that STA is much more appropriate. Fig.\ref{fig5} depicts the trajectories of parameters $K,T$, and $\tau$, and it is shown that STA can converge to the neighborhood of global solutions at 20 iterations.

Again, with the estimated parameters, an off-line PID controller is applied to this time-delay system. In this experiment,
the relative variables are given by\\
\indent \qquad \qquad \quad $K_p \in [0,1], K_i \in [0,1], K_d \in [0,1]$, $y_r=1, N=1500.$\\

\begin{figure}[htbp]
\centering
\subfloat[$K_p$]{\includegraphics[width=8cm,height=6cm]{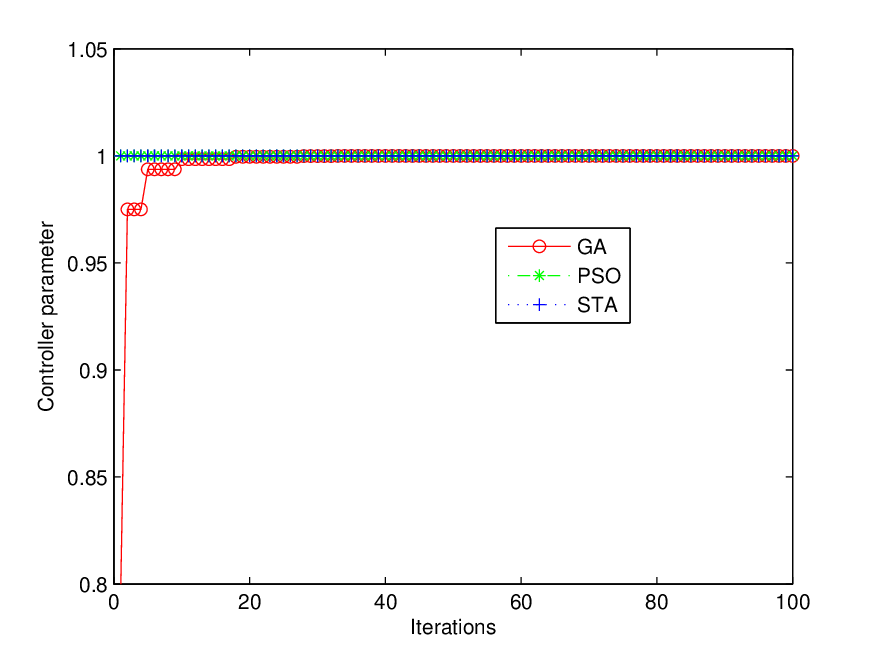}}
\subfloat[$K_i$]{\includegraphics[width=8cm,height=6cm]{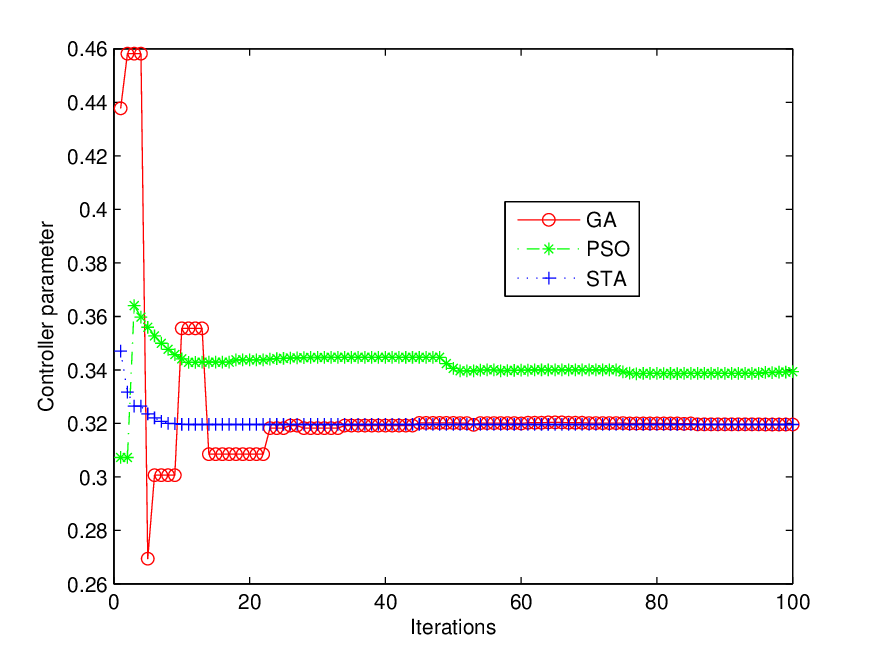}}\\
\subfloat[$K_d$]{\includegraphics[width=8cm,height=6cm]{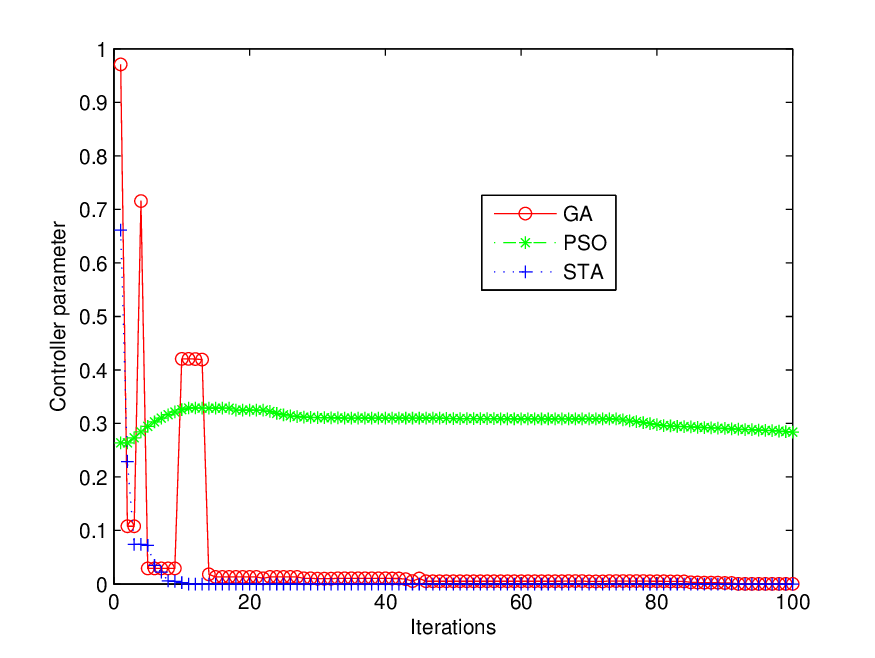}}
\subfloat[state and error]{\includegraphics[width=8cm,height=6cm]{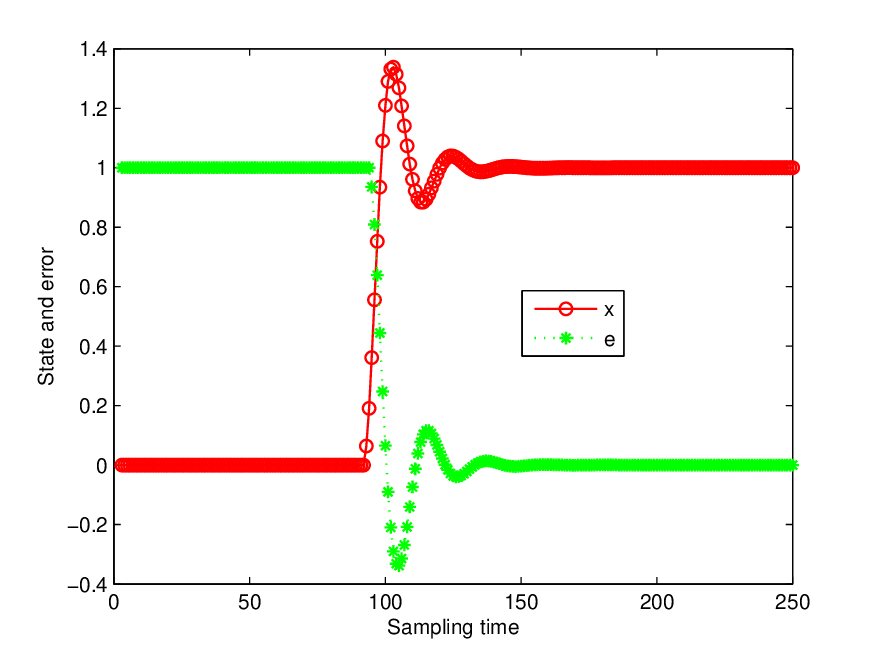}}\\
\caption{(a),(b),(c) depict the convergence of PID controller parameters $K_p,K_i$ and $K_d$ with GA, PSO and STA, respectively, and (d) shows the changes of state and error under the best parameters with STA for Example 2.}
\label{fig6}
\end{figure}

As shown in Table \ref{tab7}, the GA and STA can find a PI controller for the time-delay system, since under the PI controller,
the MSE is smaller than that of PID controller obtained by PSO. Table \ref{tab8} indicates that STA has the capacity to find the minimum MSE in a much higher probability. Fig.\ref{fig6} depicts the convergence of PID controller parameters $K_p,K_i$ and $K_d$ with GA, PSO and STA, respectively, and it is shown that STA can find the optimal parameters in a much faster way.
\begin{table}[htbp]
\begin{center}
\caption{Best parameters of PID controller for Example 2}
\label{tab7}
\begin{tabular}{{cccc}}
\hline
\toprule[1pt]
\textit{Algorithms} & $K_p$ & $K_i$ & $K_d$\\
\hline
GA  & 1.0000  &  0.3196  &   0 \\
PSO & 1.0000  &  0.3393  &   0.2837\\
STA & 1.0000  &  0.3196  &   0\\
\bottomrule[1pt]
\hline
\end{tabular}
\end{center}
\end{table}

\begin{table}[htbp]
\begin{center}
\caption{Performance comparison of PID controller for Example 2}
\label{tab8}
\begin{tabular}{{ccccc}}
\hline
\toprule[1pt]
\textit{Algorithms} & Best & Mean & Worst & St.dev\\
\hline
GA  & 6.3256E-2 & 6.3256E-2 &   6.3256E-2 &    1.0518E-7\\
PSO & 6.3258E-2 & 6.3274E-2 &   6.3279E-2 &    6.6540E-6\\
STA & 6.3256E-2 & 6.3256E-2 &   6.3256E-2 &    4.0097E-15\\
\bottomrule[1pt]
\hline
\end{tabular}
\end{center}
\end{table}

\section{Conclusion}
In this paper, a new optimization algorithm named STA is applied to solve the problems of parameter estimation and controller design for nonlinear systems. As a optimizer, STA is used to achieve the accurate model, and then it is adopted to obtain the optimal off-line PID controller. The experimental results have confirmed the validity of proposed algorithm. By comparison with GA and PSO, it is found that STA has stronger global search ability and is more stable in statistics. With regard to the convergence rate, it is also discovered that STA is much faster than its competitors. As a novel optimization method, these applications of STA show that it is a promising alternative approach for system identification and control.
\section*{Acknowledgements}
The work was supported by the National Science Found for Distinguished Young Scholars of China (Grant No. 61025015), the Foundation for Innovative Research Groups of the National Natural Science Foundation of China (Grant No. 61321003) and the China
Scholarship Council.




\begin{thebibliography}{100}

\bibitem{fei} Minrui Fei, Dajun Du and Kang Li, A fast model identification method for networked control system, Applied Mathematics and Computation, 205(2)(2008), 658--667.
\bibitem{pli} Pingkang Li, Kang Li, A recursive algorithm for nonlinear model identification, Applied Mathematics and Computation, 205(2)(2008), 511--516.
\bibitem{chen} Jing Chen, Xianling Lu, Rui Ding, Parameter identification of systems with preload nonlinearities based on the finite impulse response model and negative gradient search, Applied Mathematics and Computation, 219(5)(2012), 2498--2505.
\bibitem{jli} Junhong Li, Rui Ding, Parameter estimation methods for nonlinear systems, Applied Mathematics and Computation, 219(9)(2013), 4278--4287.
\bibitem{bi1} Alireza Alfi, Hamidreza Modares, System identification and control using adaptive particle swarm optimization, Applied Mathematical Modelling, 35(2011) 1210-1221.
\bibitem{bi2} Hamidreza Modares, Alireza Alfi, Mohammad-Bagher Naghibi Sistani, Parameter estimation of bilinear systems based on an adaptive particle swarm optimization, Engineering Applications of Artificial Intelligence, 23(2010) 1105-1111.
\bibitem{bi3} Wei-Der Chang, Nonlinear system identification and control using a real-coded genetic algorithm, Applied Mathematical Modelling, 31(2007) 541-550.
\bibitem{bi4} K.J., Astrom, B. Wittenmark, Adaptive Control, Addison-Wesley, Massachusetts, 1995.
\bibitem{bi5} ALFI Alireza, PSO with adaptive mutation and inertia weight and its application in parameter estimation of dynamic Systems, ACTA ATUOMATICA SINICA, 37(5)(2011), 541-549.
\bibitem{bi6} Goldberg D E, Genetic Algorithm in Search, Optimization and Machine Learning, NJ: Addison Wesley, 1989.
\bibitem{bi7} Holland, J.H, Adaptation in Natural and Artificial Systems, Michigan: The University of Michigan Press, 1975.
\bibitem{bi8} Eberhart, R.C. and Kennedy, J., A new optimizer using particles swarm theory, in Proceedings of Sixth International Symposium on Micro Machine and Human Science, 1995, pp.39-43.
\bibitem{bi9} Yuhui Shi, Russell C.Eberhart. Empirical study of particle swarm optimization, in Proceedings of IEEE International Congress on Evolutionary Computation, (3)(1999):591-600.
\bibitem{bi10} Murat \.{I}hsan K\"{o}m\"{u}rc\"{u}, Nedim Tutkun, \.{I}sma\.{i}l Hakk{\i} \"{O}z\"{o}l\c{c}er, Adem Akp{\i}nar, Estimation of the beach bar parameters using the genetic algorithms, Applied Mathematics and Computation, 195(2008),49-60.
\bibitem{bi11} Dionysios C. Aliprantis, Scott D. Sudhoff, Brian T. Kuhn, Genetic algorithm-based parameter identification of a hysteretic brushless exciter model, IEEE Transaction on Energy Conversion, 21(1)(2006), 148-154.
\bibitem{bi12} K. Valarmathi, D. Devaraj, T.K. Radhakrishnan, Real-coded genetic algorithm for system identification and controller tuning, Applied Mathematic Modelling, 33(2009) 3392-3401.
\bibitem{bi13} Romain Marion, Riccardo Scorretti, Nicolas Siauve, Marie-Ange Raulet and Laurent Kr\"{a}henb\"{u}hl, Identification of Jiles-Atherton model parameters using particle swarm optimization, IEEE Transaction on Magnetics, 44(6)(2008), 894-897.
\bibitem{bi14} Qi Li, Weirong Chen, Youyi Wang, Shukui Liu and Junbo Jia, Parameter identification for PEM fuel-cell
    mechanism model based on effective informed adaptive particle swarm optimization, IEEE Transations on Industrial
    Electronics, 58(6)(2011), 2410-2419.
\bibitem{bi15} R. Eberhart, Y. Shi, Comparison between genetic algorithms and particle swarm optimization, Annual Conference on Evolutionary Programming, San Diego, 1998.
\bibitem{bi16} M. Clerc, J. Kennedy, The particle swarm-explosion, stability, and convergence in a multidimensional complex space, IEEE Transactions on Evolutionary Computaion, 6(1)(2002), 58-73.
\bibitem{bi17} X.J. Zhou, C.H. Yang and W.H. Gui, Initial version of state transition algorithm, in the 2nd International Conference on Digital Manufacturing and Automation, 2011, pp.644--647.
\bibitem{bi18} X.J. Zhou, C.H. Yang and W.H. Gui, A new transformation into state transition algorithm for finding the global minimum, in the 2nd International Conference on Intelligent Control and Information Processing, 2011, pp.674--678.
\bibitem{xzhou2012} X.J. Zhou, C.H. Yang and W.H. Gui, State transition algorithm, Journal of Industrial and Management Optimization, 2012, 8(4): 1039--1056.
\bibitem{xzhou2013} X.J. Zhou, D.Y. Gao, C.H. Yang, A Comparative study of state transition algorithm with harmony search and artificial bee colony, Advances in Intelligent Systems and Computing, 212(2013), 651--659.
\bibitem{bi19} F.G. Shinskey, Process Control System: Application, Design and Tuning, McGraw-Hill, 1996.
\bibitem{bi20} T.L. Seng, M.B. Khalid, R.Yusof, Tuning of a neuro-fuzzy controller by genetic algorithm, IEEE Transactions on System, Man and Cybernetics(B), 29(1999), 226-236.


\end{thebibliography}



\end{document}